\documentclass[12pt]{article}
\usepackage[cp866]{inputenc}
\usepackage[english]{babel}
\usepackage[psamsfonts]{amssymb}
\usepackage{amsmath,amsfonts,amssymb,amscd}
\newtheorem{theorem}{Theorem}[section]

\newtheorem{lemma}[theorem]{Lemma}
\newtheorem{proposition}[theorem]{Proposition}
\newtheorem{question}[theorem]{Question}

\begin{document}

\title{Example of non-linearizable quasi-cyclic subgroup of automorphism group of polynomial algebra
\footnotetext{This research  was partially supported by
the Laboratory of Quantum Topology of Chelyabinsk State University (Russian Federation government grant 14.Z50.31.0020),
 RFBR-14-01-00014,  RFBR-13-01-00513
and Indo-Russian RFBR-13-01-92697.
}}

\author{Valeriy G. Bardakov\\
Sobolev Institute of Mathematics and Novosibirsk State University, \\
Novosibirsk 630090, Russia\\
Laboratory of Quantum Topology, Chelyabinsk State University, \\
Brat'ev Kashirinykh street 129, Chelyabinsk 454001, Russia.\\
bardakov@math.nsc.ru\\
\vspace{5pt}\\
Mikhail V. Neshchadim
\\Sobolev Institute of Mathematics and Novosibirsk State University, \\
Novosibirsk 630090, Russia\\
neshch@math.nsc.ru}

\date{}

\markboth{V.\,G.~Bardakov, M.\,V.~Neshchadim}{Example of non-linearizable quasi-cyclic subgroup}

\maketitle

\begin{abstract}
It is well known that every finite subgroup of automorphism group of polynomial algebra of rank 2 over the field of zero characteristic is conjugated with a subgroup of linear automorphisms.
We  prove that it is not true for an arbitrary torsion subgroup. We construct an example of abelian  $p$-group of automorphism of polynomial algebra of rank 2 over the field of complex numbers, which is not conjugated with a subgroup of linear automorphisms.


\end{abstract}


\subsection{Introduction}
\label{subsec1}

In the present paper we consider the polynomial algebra
$P_n = k [ x_1, x_2, \ldots, x_n ]$, and free  associative algebra
$A_n = k \langle x_1, x_2, \ldots, x_n \rangle$
 with free generators
$x_1, x_2, \ldots, x_n $ over the field $k$. We assume that these algebras contain unit elements.
We will use symbols $\mathrm{Aut} \, P_n$ and $\mathrm{Aut} \, A_n$, to denote groups of $k$-automorphisms of these algebras, i.~e. automorphisms which fix elements of $k$.
It is well known  \cite{C}, that the group $\mathrm{Aut} \, P_2$ is isomorphic to the group  $\mathrm{Aut} \, A_2$.
Let $\widetilde{P}_n = k [[ x_1, x_2, \ldots, x_n ]]$ be the algebra of the formal power series in commutative variables. The following inclusions evidently holds
$$
P_n \subseteq \widetilde{P}_n,\,\,\,\,\,
\mathrm{Aut} \, P_n \subseteq \mathrm{Aut} \, \widetilde{P}_n.
$$

In the theory of automorphisms of algebra $P_n$ the following problem is  well known (see the survey \cite{K}): Is it true that every automorphism  from
 $\mathrm{Aut} \, P_n$ of finite order is conjugated  with linear automorphism in  $\mathrm{Aut} \, P_n$?
Corresponding question for involutions was formulated in   \cite[Question 14.68]{KT}.
The case of involutions is interesting because of its connection with the Cancelation Problem
 (see for example \cite{BNS}).

For  $n = 2$ and the field $k = \mathbb{C}$  of complex numbers the answer is positive, and it follows from the result of  van der Kulk
 \cite{Ku} that the group
$\mathrm{Aut} \, P_2$ is a free product with amalgamation. For  subgroups of $\mathrm{Aut} \, P_2$ more general result is known:
every finite subgroup is conjugated with the subgroup of linear automorphisms.
 P.\,M.~Cohn \cite{Co} studied algebras over the finite field and proved that every finite subgroup  of automorphism group of
$A_2$ is conjugate with subgroup of linear automorphisms if the characteristic of the field does not divide the order of the subgroup.
As reviewer noted this result was firstly proved in  \cite{Iga}. Unfortunately, this publication is not available for wide range of readers.

On the other hand, we know examples of finite subgroups of $\mathrm{Aut} \, P_3$, which are not conjugated with subgroups of linear automorphisms (for example, dihedral group of order $8 m$ for $m \geq 3$) \cite{MMP}. H.~Kraft and G.~Schwarz \cite[p. 61]{K} have formulated the following question: is there exists a commutative  (reductive) subgroup of  $\mathrm{Aut} \, P_n$, which is not conjugated with subgroup of linear automorphisms?

M. A. Shevelin \cite{Sh2} constructed an example of infinite torsion subgroup of automorphism group of free Lie algebra of rank 3 over the complex field, which is not conjugated with subgroup of linear automorphisms.
Using this idea we prove  the following assertion.

\medskip

\begin{theorem}
Let $p$ be a prime number. There are an uncountable set of pairwise non-conjugated subgroups of the group $\mathrm{Aut} \, P_2$, which are isomorphic to the quasi-cyclic
$p$--group $\mathbb{C}_{p^{\infty}}$ and which are not conjugated with subgroups of linear automorphisms.
\end{theorem}
\medskip

We are greatfull to all the participants of the  scientific
seminar  ``Evariste Galois'' (Novosibirsk State University)
for the useful discussions.  Also, we thank anonymous reviewer for helpful notes.

\subsection{Quasi-cyclic subgroups}
\label{subsec2}

Throughout the paper we use the  field of complex numbers
 $\mathbb{C}$ as the main field. Every automorphism  $\psi \in \mathrm{Aut} \, P_n$ is completely determined by its action on the generators, therefore we identify automorphism $\psi$ with the vector $(x_1^{\psi}, x_2^{\psi}, \ldots , x_n^{\psi})$, i.~e.
$$
\psi = (x_1^{\psi}, x_2^{\psi}, \ldots , x_n^{\psi}).
$$

Let $p$ be a fixed prime number. Recall  (see, for example, \cite[1.2.5, 2.4.10]{KM}) that a
{\it quasi-cyclic  $p$--group} $\mathbb{C}_{p^{\infty}}$ is the set of all complex solutions of the equations $x^{p^n} = 1$,
$n = 1, 2, \ldots$ with natural multiplication. Group $\mathbb{C}_{p^{\infty}}$ is an infinite abelian  $p$--group, such that every proper subgroup of it is cyclic.
Let us denote
$$
G_p = \{ (\alpha x_1, \alpha x_2) \, |\, \alpha \in \mathbb{C}_{p^{\infty}} \}.
$$
It is obvious that $G_p$ is a subgroup of linear automorphisms in $\mathrm{Aut} \, P_2$ and it is isomorphic to $\mathbb{C}_{p^{\infty}}$.

We will prove that there exits other embedding of
 $\mathbb{C}_{p^{\infty}}$ into the group $\mathrm{Aut} \, P_2$.
To do it consider a formal power series
$$
\sum_{k=0}^{\infty} a_k \, w_k, \,\,\, a_k \in \mathbb{C},\,\,\, w_k = x_2^{p^k + 1}
$$
in the algebra $\widetilde{P}_2$.
Then the map $a : \widetilde{P}_2 \longrightarrow \widetilde{P}_2$, which is defined by the equality
$$
a = \left( x_1 + \sum_{k=0}^{\infty} a_k \, w_k, \, x_2  \right)
$$
is an automorphism of $\widetilde{P}_2$. If there are only a finite number of non-zero coefficients $a_k$ here, then
$a$ is an automorphism of the algebra $P_2$.

Consider the subgroup
$$
G_p^a = a^{-1} \, G_p \, a
$$
of the group $\mathrm{Aut} \, \widetilde{P}_2$. Let us show that $G_p^a$ lies in the group
$\mathrm{Aut} \, P_2$, i.~e. it is a subgroup of automorphism group of  $P_2$.
Really, let an automorphism  $\varphi = (\alpha x_1, \alpha x_2)$, $\alpha \in \mathbb{C}_{p^{\infty}}$ belongs to the group $G_p$.
Conjugating it by automorphism  $a$ we have an equality
$$
a^{-1} \varphi a = \left( \alpha x_1 + \alpha \sum_{k=0}^{\infty} a_k \, (1 - \alpha^{p^k}) \, w_k, \,
\alpha x_2
\right).
$$
Here and further we assume that automorphisms act from the left to the right.
Since $\alpha \in \mathbb{C}_{p^{\infty}}$ we see that for sufficiently large number $k_0$ all the coefficients
 $1 - \alpha^{p^k}$, $k = k_0, k_0+1, \ldots$ are equal to zero and hence  $a^{-1} \varphi a \in \mathrm{Aut} \, P_2$.

By construction, the group  $G_p^a$ is conjugated with the group $G_p$ in the automorphism group
$\mathrm{Aut} \, \widetilde{P}_2$ of the formal power series  $\widetilde{P}_2$.
Let us find out the question about conjugation of  $G_p^a$ with some subgroup of linear automorphisms in the group
$\mathrm{Aut} \, P_2$.

\medskip

\begin{proposition}
If an automorphism
$$
a = \left( x_1 + \sum_{k=0}^{\infty} a_k \, w_k, \, x_2  \right)
$$
is defined by infinite number of non-zero coefficients  $a_k$, then $G_p^a$  is not conjugated with any subgroup of linear automorphisms.
\end{proposition}

\medskip

{\large Proof.}
Let $H$ be a subgroup of linear automorphisms, which is conjugated with $G_p^a$, i.~e.
$$
\theta^{-1} \, G_p^a \, \theta = H,
$$
for some $\theta \in \mathrm{Aut} \, P_2$. With no loss of generality we assume that $H$ belongs to the group of diagonal automorphisms, therefore the last equality means that for every automorphism
$\varphi = (\alpha x_1, \alpha x_2)$
there exists an automorphism  $h = (\alpha_1 x_1, \alpha_2 x_2) \in H$,
$\alpha_1, \alpha_2 \in \mathbb{C}_{p^{\infty}}$, such that the following equality holds
\begin{equation}
\varphi^a \theta = \theta h.
\label{1}
\end{equation}
Suppose that
$
\theta = (f_1, f_2),
$
where the polynomials  $f_1 = f_1(x_1, x_2)$ and  $f_2 = f_2(x_1, x_2)$ have degree $m_1$ an $m_2$ correspondingly.
Since these polynomials generate the algebra  $P_2$, then
$m_1 > 0$ and  $m_2 > 0$.
It is not difficult to see that the equality   (\ref{1}) is equivalent to the system
$$
 \alpha f_1 + \alpha \sum_{k=0}^{\infty} a_k \, (1 - \alpha^{p^k}) \, f_2^{p^k+1}    =
f_1(\alpha_1 x_1, \alpha_2 x_2),  \,\,\,\,
\alpha f_2(x_1, x_2) = f_2(\alpha_1 x_1, \alpha_2 x_2).
$$
If we write the first equality in the following form
$$
 \alpha \sum_{k=0}^{\infty} a_k \, (1 - \alpha^{p^k}) \, f_2^{p^k+1} =
 f_1(\alpha_1 x_1, \alpha_2 x_2) -  \alpha f_1,
$$
then we see that the right hand side of it has the degree, which is less or equal to $m_1$. But taking different
$\alpha$, we can get an arbitrary large degree of the left hand side of this equality. Hence, the equality (\ref{1})
is impossible. The proposition is proved.

\subsection{On conjugation of quasi-cyclic subgroups}
\label{subsec3}

In this section we will prove the following proposition.

\begin{proposition}
There exist an uncountable set of subgroups  $G_p^a$, which are pairwise non-conjugated in $\mathrm{Aut} \, P_2$.
\end{proposition}

Let
$$
a = \left( x_1 + \sum_{k=0}^{\infty} a_k \, w_k, x_2  \right),\,\,\,
b = \left( x_1 + \sum_{k=0}^{\infty} b_k \, w_k, x_2  \right),
$$
$a_k,\,b_k \in \mathbb{C}$, $k = 0, 1, 2, \ldots$
be two automorphisms of  $\widetilde{P}_2$, which are not automorphisms of  $P_2$, i.~e.
each of them contains an infinite number of non-zero coefficients.
Let us find out are the subgroups $G_p^a$ and  $G_p^b$ conjugated in  $\mathrm{Aut} \, P_2$.

Let $\theta$ be such an automorphism of $P_2$, that the following equality holds
$$
\theta^{-1} (a^{-1} \varphi a) \theta = b^{-1} \varphi b \,\, \mbox{for all}\,\, \varphi \in G_p.
$$
This equality is evidently equivalent to the following equality
\begin{equation}
 (a^{-1} \varphi a) \theta = \theta (b^{-1} \varphi b).
 \label{2}
\end{equation}
Suppose that the automorphism  $\theta$ acts on the generators of $P_2$ by the following way
$$
x_1^{\theta} = f_1 (x_1, x_2),\,\,\, x_2^{\theta} = f_2 (x_1, x_2).
$$
Then the equality (\ref{2}) is equivalent to the system
$$
 \left( \alpha x_1 + \alpha \sum_{k=0}^{\infty} a_k \, (1 - \alpha^{p^k}) \, w_k  \right)^{\theta}  =
f_1^{b^{-1} \varphi b},\,\,\,\,
(\alpha x_2)^{\theta} = f_2^{b^{-1} \varphi b},
$$
which can be rewritten by the following way

\begin{equation}\label{3}
 \left\{
\begin{array}{l}
 \alpha f_1 + \alpha \sum_{k=0}^{\infty} a_k \, (1 - \alpha^{p^k}) \, f_2^{p^k+1}    =
f_1 \left( \alpha x_1 + \alpha \sum_{k=0}^{\infty} b_k \, (1 - \alpha^{p^k}) \, w_k, \,  \alpha x_2 \right),  \\
\\
\alpha f_2 = f_2 \left(\alpha x_1 + \alpha \sum_{k=0}^{\infty} b_k \, (1 - \alpha^{p^k}) \, w_k, \, \alpha x_2 \right),        \\
\end{array}
\right.
\end{equation}
By the choice of $\alpha$,  the following equalities holds $1 - \alpha^{p^k} = 0$ for enough large $k$. Hence, the system of equations
 (\ref{3}) is the system of equalities in the algebra $P_2$.

If the polynomial $f_2$ contains the variable $x_1$, then it is decomposed  by the powers of  $x_1$
$$
f_2 = \sum_{j=0}^s c_j(x_2) \, x_1^j,\,\,\,\, c_j(x_2) \in \mathbb{C}[x_2].
$$
Then
$$
f_2 \left(\alpha x_1 + \alpha \sum_{k=0}^{\infty} b_k \,  (1 - \alpha^{p^k}) \, w_k, \, \alpha x_2 \right) =
\sum_{j=0}^s c_j(\alpha x_2) \left(\alpha x_1 + \alpha \sum_{k=0}^{\infty} b_k \,
(1 - \alpha^{p^k}) \, w_k \right)^j
$$
and since  $\alpha \in \mathbb{C}_{p^{\infty}}$ is arbitrary, then the degree of the right hand side with respect to $x_2$ can be arbitrarily large. Hence, $f_2$ does not contain  $x_1$,
i.~e. $f_2 = f_2(x_2) \in \mathbb{C}[x_2]$.
Since $\theta$ is an automorphism, then
$$
f_2 = \beta x_2 + \beta_0,\,\,\, \beta \in \mathbb{C}^*,\,\, \beta_0 \in \mathbb{C},
$$
and the first equality of  (\ref{3}) has  the following form
\begin{equation}
\begin{array}{l}
 \alpha f_1 (x_1, x_2) + \alpha \sum_{k=0}^{\infty} a_k \,(1 - \alpha^{p^k}) \, (\beta x_2 + \beta_0)^{p^k+1} =\\ \\
=f_1 \left( \alpha x_1 + \alpha \sum_{k=0}^{\infty} b_k \, (1 - \alpha^{p^k}) \, w_k, \, \alpha x_2 \right).
\end{array}
\label{4}
\end{equation}
Since $\theta = (f_1, \beta x_2  + \beta_0)$ we see that
$$
f_1 = \gamma x_1 + g(x_2),\,\,\,\, \gamma \in \mathbb{C}^*,\,\,\, g \in \mathbb{C}[x_2].
$$
Then the equality  (\ref{4}) has the form
$$
\begin{array}{l}
 \alpha (\gamma x_1 + g(x_2)) +
 \alpha \sum_{k=0}^{\infty} a_k \, (1 - \alpha^{p^k}) \, ( \beta x_2  + \beta_0)^{p^k+1}    =\\ \\
 =\gamma \left( \alpha x_1 + \alpha \sum_{k=0}^{\infty} b_k \,
  (1 - \alpha^{p^k})  \, w_k  \right) + g(\alpha x_2).
 \end{array}
$$
or equivalently
$$
 \alpha g(x_2) - g(\alpha x_2) + \alpha \sum_{k=0}^{\infty}   (1 - \alpha^{p^k}) \, [a_k (\beta x_2 + \beta_0)^{p^k+1} - b_k x_2^{p^k+1}]
  = 0.
$$
Since the degree of  $g(x_2)$ is bounded and there exist an infinite number of non-zero coefficients  $a_k$ and  $b_k$, then the following equalities hold for some number
$k_0$ which is large enough
$$
a_k \beta^{p^k+1} = b_k,\,\,\, k \geq k_0.
$$
If it is not true, then there is no automorphism  $\theta$ of $P_2$ such that
$$
\theta^{-1} (a^{-1} \varphi a) \theta = b^{-1} \varphi b.
$$
Thus we have proved that if the sequences  $\{ a_k \}$ and  $\{ b_k \}$ satisfy the  inequality
 $$
a_k \beta^{p^k+1} \not = b_k
$$
 for an infinite set of indexes $k$ and for an arbitrary $\beta\in \mathbb{C}^{*}$, then the subgroups  $G_p^a$ and  $G_p^b$ are not conjugated in  $\mathrm{Aut} \, P_2$.

To finish the proof of the proposition we have to prove that there exist an uncountable set of such sequences.

Denote by
$$
\Omega = \{ \lambda = (\lambda_1, \lambda_2, \lambda_3, \ldots )\, |\, \lambda_i \in \{ 0, 1 \} \}
$$
the set of all infinite binary sequences. The set $\Omega$ evidently contains an uncountable set of different sequences. Furthermore, the following lemma is true

\medskip

\begin{lemma}
There exists an uncountable subset $\Omega_0 \subset \Omega$, which posses the following property: for any two sequences $\lambda, \mu \in \Omega_0$ inequality
$\lambda_k \not= \mu_k$ holds for an infinite number of indexes  $k \in \mathbb{N}$.
\end{lemma}

{\large Proof.}
If we define the following relation on the set $\Omega$
$$
\lambda \sim \mu  \Leftrightarrow \lambda \,\,
 \mbox{and}\,\, \mu  \,\, \mbox{different only on finite number of places},
$$
then it is not difficult to see that this relation is equivalence relation and hence, divides the set  $\Omega$ onto equivalence classes. Moreover, every equivalence class contains only countable number of elements and then the number of equivalence classes
 $\Omega /{\sim}$ is uncountable. The transversal set with respect to relation $\sim$ is the necessary subset $\Omega_0$.

It follows from this lemma, that there is an uncountable set of subgroups
$G_p^a$, which are pairwise  non-conjugated in  $\mathrm{Aut} \, P_2$ . It completes the proof of Proposition 2.

\medskip

Now the main theorem follows from the propositions 1 and 2.

\medskip

We have found a sufficient condition when subgroups  $G_p^a$ and $G_p^b$ are not conjugated in $\mathrm{Aut} \, P_2$.

\medskip

\begin{question} Find   necessary and sufficient conditions when subgroups $G_p^a$ and $G_p^b$ are conjugated in
$\mathrm{Aut} \, P_2$.
\end{question}


\begin{thebibliography}{HD}


\bibitem{C} A.\,J.~Czerniakiewicz
\textit{Automorphisms of a free associative algebra of rank 2. Part I}, Trans. Amer. Math. Soc.  160 (1977), 393--401.


\bibitem{K}
H.~Kraft, G.~Schwarz
\textit{Finite automorphisms of affine $N$-space. Automorphisms of affine spaces}
Procceedings of a Conference held in Curacao (Netherlands Antilles), July 4--8, 1994
 Kluwer Academic Publishers, 1995,  55--66.

\bibitem{KT}
The Kourovka Notebook, Unsolved Problems in Group Theory, 18th ed., Sobolev Institute of Mathematics,
Novosibirsk, 2014.


\bibitem{BNS}
 V.\,G.~Bardakov, M.\,V.~Neshchadim, Yu.\,V.~Sosnovsky
\textit{Groups of triangular automorphisms of
a free associative algebra and a polynomial algebra}. J. Algebra
362 (2012), 201--220.



\bibitem{Ku}
 W.~van der Kulk
\textit{On polynomial rings in two variables}.
Nieuw Archief voor Wiskunde,  3:1, (1953), 33--41.



\bibitem{MMP}
 M.~Masuda, L.~Moser-Jauslin, T.~Petrie
\textit{Equivariant algebraic vector bundles over
representations of reductive groups: applications}.
Proc. Nat. Acad. Sci. U.S.A., 88:20 (1991), 9065--9066.



\bibitem{Co}
 P.\,M.~Cohn
\textit{The  automorphism group of the free algebra of rank two}. Serdica Math. J., 28 (2002), 255--266.

\bibitem{Iga}
T.~Igarashi
\textit{Finite subgroups of the  automorphism group of the affine plain}.
 Thesis, Osaka University,  1977.



\bibitem{Sh2}
M.\,A.~Shevelin
\textit{Subgroups of automorphism group of free Lie algebra of rank 3}. (Russian) Sibirsk. Mat. Zh. 55 (2014).


\bibitem{KM}
M.\,I.~Kargapolov,  Yu.\,I.~Merzljakov
\textit{Fundamentals of the theory of groups.} Translated from the second Russian edition by Robert G. Burns. Graduate Texts in Mathematics, 62. Springer-Verlag, New York-Berlin, 1979. xvii+203 pp. ISBN: 0-387-90396-8 20-01.

\end{thebibliography}
\end{document}